\documentclass[12pt]{amsart}
\usepackage{amssymb}
\usepackage{amsfonts,amsthm,amsmath,amscd}
\begin{document}
\hfuzz=35pt
\title[]{\it{Correction: Stable homotopy classification of} $BG_{p}^{\wedge}$}
\author{John Martino}
\address{Department of Mathematics \\ Western Michigan University
\\ Kalamazoo, MI 49008}
\email{john.martino@wmich.edu}
\author{Stewart Priddy }
\address{Department of Mathematics\\ Northwestern University\\ Evanston, IL 60208}
\email{priddy@math.northwestern.edu}
\date{\today}
\bigskip
\bigskip
\bigskip

\maketitle

\newtheorem{Thm}{Theorem}
\newtheorem{Lem}[Thm]{Lemma}
\newtheorem{Cor}[Thm]{Corollary}
\newtheorem{Prop}[Thm]{Proposition}
\newtheorem{Ex}[Thm]{Example}

\medskip
In this note we correct the proof of Theorem 1.1 of our
paper \cite{mp}, given as Theorem 1 below. Let $G$, $G'$ be finite groups and $p$ be a prime number. The goal
is to give necessary and sufficient algebraic conditions on the $p$-subgroups of $G$
and $G'$ which determine if their $p$-completions $BG^{\wedge}_p$ and $BG'^{\wedge}_p$
are stably homotopy equivalent.

\medskip

\begin{Thm}
For finite groups  $G$, $G'$ the following are equivalent:

\begin{enumerate}
\item {$BG^{\wedge}_p$ and $BG'^{\wedge}_p$ are stably homotopy equivalent.}

\item For every $p$--group $Q$,
$$ \mathbf F_p Rep(Q,G) \approx \mathbf F_p Rep (Q,G') $$
as $Out(Q)$ modules.  $Rep(Q,G) = Hom(Q,G)/G$ with G acting by conjugation.

\item For every $p$--group $Q$,
$$ \mathbf F_p Inj(Q,G) \approx \mathbf F_p Inj(Q,G') $$
as $Out(Q)$ modules.  $Inj(Q,G) \subseteq Rep(Q,G)$ consists of conjugacy
classes of injective homomorphisms.

\end{enumerate}
\end{Thm}

Our notation follows that of \cite{mp}.
Throughout all groups are finite and hereafter all group maps are monomorphisms.
\medskip

The key step in the proof is the next Proposition.
For $Q$ a $p$-group we let
$Cen(Q,G) \subseteq Inj(Q,G)$ be the equivalence classes of injections
$\alpha :Q\longrightarrow G$ such that $C_{G}(Im\alpha )/Z(Im\alpha )$
is a $p'$-group.

\begin{Prop} If $\mathbf{F}_p Inj(Q,G) \approx \mathbf{F}_p Inj(Q,G')$ as $Out(Q)$ modules for all
$p$-groups $Q$ then  $\mathbf{F}_pCen(Q,G) \approx \mathbf{F}_pCen(Q,G')$ as $Out(Q)$ modules for all $Q$.
\end{Prop}

\medskip

Let $nCen(Q,G)$ be the complement of $Cen(Q,G)$ in $Inj(Q,G)$.
To study these sets we shall consider subgroups $H\leq G$
such that $H\approx Q$ and $p$ divides $|C_G (H)/Z(H)|$.
Let $\widetilde{H}$ be a Sylow
$p$-subgroup of $H\cdot C_G (H)$. By definition $|H| < |\widetilde{H}|$.
An equivalence $s: \widetilde{H}_1
\overset{\approx}\longrightarrow \widetilde{H}_2$
between such groups will be an isomorphism such that $s(H_1) = H_2$.
Let $\{\widetilde{Q}_{j} \}$ be a set
of representatives for equivalence classes of such subgroups.

The proof of Proposition 2 depends on a close analysis of $Cen(Q,G)$. Define
 $$ Cen_{Q_j}(\widetilde{Q}_j, \widetilde{Q}_k, G) = \{ \beta: \widetilde{Q}_j \to G \  | \  
[\beta] \in Cen(\widetilde{Q}_j, G), \widetilde{\beta(Q_j)} \sim \widetilde{Q}_k \}/G. $$

The following lemma describes some elementary properties of $ Cen_{Q_j}(\widetilde{Q}_j, \widetilde{Q}_k, G)$.
Proofs are given at the end of the note.

\begin{Lem} Let $\widetilde{Q}_m$ denote a subgroup of maximum order in $\{\widetilde{Q}_j\}$.
\begin{enumerate}
\item $\{ \beta: \widetilde{Q}_j \to G \ | \  \widetilde{\beta(Q_j)} \sim \widetilde{Q}_j \}/G 
\subseteq Cen(\widetilde{Q}_j , G).$
\item $Cen(\widetilde{Q}_m , G) = Inj(\widetilde{Q}_m,G).$
\item $Cen_{Q_m}(\widetilde{Q}_m, \widetilde{Q}_m, G) = Cen(\widetilde{Q}_m , G).$
\item If $|\widetilde{Q}_j| = |\widetilde{Q}_k|$ and  $Cen_{Q_j}(\widetilde{Q}_j, \widetilde{Q}_k, G) \neq \phi$
then $\widetilde{Q}_j \sim \widetilde{Q}_k$.
\item  $Cen(\widetilde{Q}_j, G) = Cen_{Q_j}(\widetilde{Q}_j, \widetilde{Q}_j, G) \coprod_k Cen_{Q_j}(\widetilde{Q}_j, \widetilde{Q}_k, G)$
where $k$ ranges over $\{k : |\widetilde{Q}_k| > |\widetilde{Q}_j|\}$.

\end{enumerate}
\end{Lem}
\medskip

Define a set of monomorphisms
$$R(Q',\widetilde{Q}) = \{\alpha:Q' \to \widetilde{Q} \ |\   \alpha(Q') = Q \}.$$
and a group of automorphisms
$$ Aut(\widetilde{Q}|Q) = \{\alpha: \widetilde{Q} \to \widetilde{Q} \ \ | \ \alpha(Q)=Q\}.$$
Observe that $R(Q',\widetilde{Q})$ and $Cen_{Q}(\widetilde{Q}, \widetilde{Q}_k, G)$ are $Aut(\widetilde{Q} |Q)$ sets since the
image of $Q$ is invariant. The following Lemma is a corrected version of Lemma 4.6 [MP] which contained an error
pointed out to us by Kari Ragnarsson \cite{r}.

\begin{Lem} There is an isomorphism of $Out(Q)$ sets
$$ \coprod_k \psi: R(Q,\widetilde{Q}_k) {\times}_{Aut(\widetilde{Q}_k|Q_k)}Cen_{Q_k}(\widetilde{Q}_k, \widetilde{Q}_k, G)
\to nCen(Q,G)$$
natural in $G$ given by composition, i.e. $\psi(\gamma \times [\delta]) = [\delta\gamma]$.
\end{Lem}

\bigskip

Define a set of monomorphisms
$$R_{{Q}_i}(\widetilde{Q}_i,\widetilde{Q}_k) = \{\alpha:\widetilde{Q}_i \to \widetilde{Q}_k \ |\   \alpha(Q_i) = Q_k \}.$$
Arguing as in Lemma 4 we have

\begin{Lem} There is an isomorphism of $Aut(\widetilde{Q}_i|Q_i)$ sets
$$ \psi: R_{Q_i}(\widetilde{Q}_i,\widetilde{Q}_k)\times_{Aut(\widetilde{Q}_k|{Q}_k)} Cen_{{Q}_k}(\widetilde{Q}_k,\widetilde{Q}_k, G)
\to Cen_{{Q}_i}(\widetilde{Q}_i, \widetilde{Q}_k,G).$$
natural in $G$ where $\psi(\gamma, [\delta])= [\delta \gamma]$.
\end{Lem}
\medskip

\begin{Prop} If  $\mathbf{F}_p Inj(Q,G) \approx \mathbf{F}_p Inj(Q,G')$ as $Out(Q)$ modules for all
$p$-groups $Q$ then $\mathbf{F}_pCen(\widetilde{Q}_j, G) \approx 
\mathbf{F}_pCen(\widetilde{Q}_j, G')$ as $Out(\widetilde{Q}_j)$ modules for all $j$ and
$$ \mathbf{F}_pCen_{Q_j}(\widetilde{Q}_j, \widetilde{Q}_j,G) \approx \mathbf{F}_pCen_{Q_i}(\widetilde{Q}_j,\widetilde{Q}_j, G')$$ as $Aut(\widetilde{Q_j}|Q_j)$ modules for all $j$.
\end{Prop}

\bigskip
{\noindent{\bf Proofs.}}

\medskip
{\it Proof of Lemma 3:} (1) Let $\beta \in \{ \beta: \widetilde{Q}_j \to G \ | \  \widetilde{\beta(Q_j)} \sim 
\widetilde{Q}_j \}$. By definition we must show
$C_G(\beta(\widetilde{Q}_j))/Z(\beta(\widetilde{Q}_j))$ is a $p'$-group. Let $x \in  C_G(\beta(\widetilde{Q}_j))
\subseteq C_G(\beta(Q_j))$ have order a power of $p$. Since $\widetilde{\beta(Q_j)} \sim 
\widetilde{Q}_j$, $|\widetilde{Q}_j| = |\beta(\widetilde{Q}_j)| \leq 
|\widetilde{\beta(Q_j)}| = |\widetilde{Q}_j|$. Therefore $\beta(\widetilde{Q}_j) = \widetilde{\beta(Q_j)}$ and so 
$\beta(\widetilde{Q_j})$ is a Sylow $p$-subgroup of $C_G(\beta(Q_j))\cdot \beta(Q_j)$.
However
$\langle x, \beta(\widetilde{Q}_j)\rangle$ is a $p$-subgroup of $C_G(\beta(Q_i))\cdot Q_i$ which means
$x \in \beta(\widetilde{Q}_j)$ and therefore $x \in Z(\beta(\widetilde{Q}_j))$.

(2) Since $Cen(\widetilde{Q}_m , G) \subseteq Inj(\widetilde{Q}_m,G)$ we need only show the other inclusion.
Let $[\beta] \in Inj(\widetilde{Q}_m,G)$. Thus we must show
$C_G(\beta(\widetilde{Q}_m))/Z(\beta(\widetilde{Q}_m))$ is a $p'$-group. Using maximality this argument is
analogous to that of (1).

(3) By definition $Cen_{Q_m}(\widetilde{Q}_m, \widetilde{Q}_m ,G) \subseteq Cen(\widetilde{Q}_m, G)$.
Let $[\beta] \in Cen(\widetilde{Q}_m, G)$. By maximality
we have equality $|\widetilde{Q}_m| = |\beta(\widetilde{Q}_m)| = |\widetilde{\beta(Q_m)}| = |\widetilde{Q}_m|$.
Since $\widetilde{\beta(Q_m)}$ is a Sylow $p$-subgroup of $C_G(\beta(Q_m))\cdot \beta(Q_m)$
there is an isomorphism
$$ 
c_g\beta: \widetilde{Q}_m \overset{\beta\ \approx}\longrightarrow \beta(\widetilde{Q}_m)
\overset{c_g \ \approx}\longrightarrow \widetilde{\beta(Q_m)}
$$
where $g \in C_G(\beta(Q_m))\cdot \beta(Q_m)$. Hence $\widetilde{Q}_m \sim \widetilde{\beta(Q_m)}$
and so $[\beta] \in Cen_{Q_m}(\widetilde{Q}_m, \widetilde{Q}_m, G)$.

(4) If $[\beta] \in Cen_{Q_j}(\widetilde{Q}_j, \widetilde{Q}_k, G)$ then $|\widetilde{Q}_j| = \beta(\widetilde{Q}_j)| \leq
|\widetilde{\beta(Q_j)}| = |\widetilde{Q}_k|$. Thus $\beta:\widetilde{Q}_j \overset{\approx}\longrightarrow 
\widetilde{\beta(Q_j)}$ and it follows that $\widetilde{Q}_j \sim \widetilde{Q}_k$.

(5)  $[\beta]\in Cen(\widetilde{Q}_j,G)$ iff $\beta:\widetilde{Q}_j \to G$ is a monomorphism and
$C_G(\beta(\widetilde{Q}_j))/Z(\beta(\widetilde{Q}_j))$ is a $p'$-group. Now $\widetilde{\beta({Q}_j)} \sim
\widetilde{Q}_k$ for some unique $k$. Hence $[\beta] \in Cen_{Q_j}(\widetilde{Q}_j, \widetilde{Q}_k, G)$.
The restriction on the range of $k$ follows from (4). The other inclusion is clear. 
\qed
\bigskip

{\it Proof of Lemma 4: } 
If $R(Q,\widetilde{Q}_k) \neq \phi$ then $|Q| = |Q_k| < |\widetilde{Q}_k|$.
Therefore $\phi(\gamma_k \times [\delta_l]) \in nCen(Q,G)$.
It is easily checked that $\psi$ is well-defined. 
\smallskip

\noindent{Injectivity:} 
Suppose
$\psi(\gamma_k \times [\delta_k]) = \psi(\gamma_l \times [\delta_l])$. Then for some $g\in G$
$$
 \delta_{l}(Q_l)= im(\delta_l \gamma_l) = c_{g}im(\delta_k \gamma_k) = c_{g}\delta_k(Q_k)
$$
Hence
$c_g[C_G(\delta_k(Q_l))\cdot \delta_k(Q_k)] = C_G(\delta_l(Q_l))\cdot \delta_l(Q_l)$
and similarly for their Sylow $p$-subgroups. Hence  $\widetilde{\delta_k(Q_k)}= \widetilde{\delta_l(Q_l)}$
up to  conjugation sending $\delta_k(Q_k)$ to $\delta_l(Q_l)$. Thus $\widetilde{Q_k} \sim \widetilde{Q_l}$
which implies  $\widetilde{Q_k} = \widetilde{Q_l}$ and $Q_k = Q_l$. Therefore we have a commutative diagram
$$\CD
Q @>{\gamma}_k >> \widetilde{Q}_k @>{\delta}_k >> G \\
@VidVV       @.         @Vc_gVV \\
Q @>{\gamma}_l >> \widetilde{Q}_k @>{\delta}_l >> G
\endCD
$$
i.e.,
${\gamma}_l = {{\delta}_l}^{-1} c_g {\delta}_k {\gamma}_k$. Thus
$$
\gamma_l \times [\delta_l] =  {{\delta}_l}^{-1} c_g {\delta}_k {\gamma}_k \times [\delta_l] = 
\gamma_k \times [\delta_l({\delta_l}^{-1}c_g {\delta}_k)] = \gamma_k \times [\delta_k] 
$$
since ${{\delta}_l}^{-1} c_g {\delta}_k \in Aut(\widetilde{Q}_k | Q_k)$.
\smallskip

\noindent Surjectivity: If $[f] \in nCen(Q,G)$ then there is the unique factorization
$$ 
f: Q\overset{\pi}\longrightarrow \widetilde{f(Q))}\overset{i} 
\hookrightarrow G 
$$
where $\pi$ is a monomorphism onto $f(Q)$ and $i$ is inclusion. Choose an equivalence $u: \widetilde{f(Q)} \to
\widetilde{Q}_j$ for some $j$ such that $u(f(Q)) = Q_j$ and set $\gamma = u\pi$, $\delta = iu^{-1}$.
By Lemma 3 (1),
$[\delta] \in Cen_{Q_j}(\widetilde{Q}_j, \widetilde{Q}_j, G)$. Since $\gamma \in R(Q,\widetilde{Q}_j)$ we have
$\psi(\gamma , [\delta]) = [\delta \gamma] = [f]$ as required. \qed
\bigskip

{ \it Proof of Lemma 5:} Fix $Q = Q_i$. First we note that $\widetilde{\delta\gamma(Q)} = 
\widetilde{\delta(Q_k)}\sim \widetilde{Q}_k$
Therefore $[\delta\gamma] \in Cen_{Q}(\widetilde{Q}, \widetilde{Q}_k, G)$.
It is now easy to check that $\psi$ is well-defined. Injectivity: Fix $Q = Q_i$ and suppose
$\psi(\gamma_1 \times [\delta_1]) = \psi(\gamma_2 \times [\delta_2])$. Then
$$ \delta_1\gamma_1 = c_g\delta_2\gamma_2 $$
for some $g \in G$. Hence $\gamma_2 = {\delta_2}^{-1}{c_g}^{-1} \delta_1 \gamma_1$. Thus
$$ \gamma_2 \times[\delta_2] = {\delta_2}^{-1} {c_g}^{-1}\delta_1\gamma_1 \times [\delta_2] = 
\gamma_1 \times[\delta_2({\delta_2}^{-1}{c_g}^{-1}\delta_1)] = \gamma_1 \times [\delta_1].$$
Surjectivity: Given $[f] \in Cen_{Q}(\widetilde{Q}, \widetilde{Q}_k, G)$, there is a factorization
$$  f:\widetilde{Q} \overset{\pi} \longrightarrow {f(\widetilde{Q})} \overset{i}\hookrightarrow G $$
and an equivalence $u:\widetilde{f(Q)} \to \widetilde{Q}_k$
where $\pi$ is a isomorphism onto $f(\widetilde{Q})$, $i$ is inclusion, and $u(f(Q)) = Q_k$.
Since $f(\widetilde{Q}) \leq C_G(f({Q}))\cdot f({Q})$ we have $c_g(f(\widetilde{Q})) \leq \widetilde{f(Q)}$
for some $g \in C_G(f(Q))\cdot f(Q)$.
Let $\gamma = uc_g \pi, \delta = i{c_g}^{-1}u^{-1}$. Then $\gamma \in R_Q(\widetilde{Q}, \widetilde{Q}_k),
[\delta] \in Cen_{Q_k}(\widetilde{Q}_k, \widetilde{Q}_k, G)$ and $\psi( \gamma, [\delta]) = [f]$. \qed
\medskip

{\it Proof of Proposition 6:} By Lemma 3 (2) and downward induction on the order of $\widetilde{Q}_i$ we may assume
$$\mathbf{F}_pCen(\widetilde{Q}_i,G) \approx \mathbf{F}_pCen(\widetilde{Q}_i,G')   \eqno(1)$$
as $Out(\widetilde{Q}_i)$ modules for all $\widetilde{Q}_i$ such that $|\widetilde{Q}_i| > |\widetilde{Q}_j|$
for some $j$.
Claim:
$$\mathbf{F}_pCen_{Q_i}(\widetilde{Q}_i, \widetilde{Q}_i,G) \approx \mathbf{F}_pCen_{Q_i}(\widetilde{Q}_i,\widetilde{Q}_i, G')$$
as $Aut(\widetilde{Q}_i|Q_i)$ modules for all $\widetilde{Q}_i$ such that $|\widetilde{Q}_i| > |\widetilde{Q}_j|$.
By Lemma 3 (3) we may assume by downward induction that
$$\mathbf{F}_pCen_{Q_k}(\widetilde{Q}_k, \widetilde{Q}_k,G) \approx \mathbf{F}_pCen_{Q_k}(\widetilde{Q}_k,\widetilde{Q}_k, G')$$ as $Aut(\widetilde{Q}_k|Q_k)$ modules
for all $\widetilde{Q}_k$ such that $|\widetilde{Q}_k | > |\widetilde{Q}_i|$ for some $i$ such that
$|\widetilde{Q}_i| > |\widetilde{Q}_j|$ .
Therefore by Lemma 5,
$$ \oplus_{k,k\neq i} \mathbf{F}_pCen_{{Q}_i}(\widetilde{Q}_i,\widetilde{Q}_k, G) \approx \oplus_{k, k\neq i}
\mathbf{F}_pCen_{{Q}_i}(\widetilde{Q}_i,\widetilde{Q}_k, G')$$
as $Aut(\widetilde{Q}_i|Q_i))$-modules.
By equation (1), Lemma 3 (5),
$$ \oplus_k \mathbf{F}_pCen_{{Q}_i}(\widetilde{Q}_i,\widetilde{Q}_k, G) \approx \oplus_k
\mathbf{F}_pCen_{{Q}_i}(\widetilde{Q}_i,\widetilde{Q}_k, G')$$
as $Aut(\widetilde{Q}_i|Q_i)$-modules.
By the Krull-Schmidt Theorem we have cancellation, hence
$$\mathbf{F}_pCen_{Q_i}(\widetilde{Q}_i, \widetilde{Q}_i, G) \approx 
\mathbf{F}_pCen_{Q_i}(\widetilde{Q}_i, \widetilde{Q}_i, G')$$
which proves the claim by induction.

Setting $Q = \widetilde{Q}_j$ in Lemma 4 we conclude
$$\mathbf{F}_pnCen(\widetilde{Q}_j, G) \approx \mathbf{F}_pnCen(\widetilde{Q}_j, G').$$
as $Out(\widetilde{Q}_j)$ modules.
By definition  $$  \mathbf{F}_pInj(\widetilde{Q}_j,G) = \mathbf{F}_pCen(\widetilde{Q}_j,G) \oplus 
\mathbf{F}_pnCen(\widetilde{Q}_j,G)$$ and we
have cancellation.  It follows that
$$\mathbf{F}_pCen(\widetilde{Q}_j,G) \approx \mathbf{F}_pCen(\widetilde{Q}_j,G') $$ for all $j$.
Repeating the argument of the claim we also have
$$ \mathbf{F}_pCen_{Q_j}(\widetilde{Q}_j, \widetilde{Q}_j,G) \approx \mathbf{F}_pCen_{Q_j}(\widetilde{Q}_j,\widetilde{Q}_j, G')$$
for all $j$.
\qed
\bigskip

{\it Proof of Proposition 2:} By Prop 6 and Lemma 4,
$$\mathbf{F}_pnCen({Q}, G) \approx 
\mathbf{F}_pnCen({Q}, G').$$

The result follows by the decomposition
$$  \mathbf{F}_pInj(Q,G) = \mathbf{F}_pCen(Q,G) \oplus \mathbf{F}_pnCen(Q,G)$$
and
cancellation.   \qed

\bigskip
{\it Proof of Theorem 1:} With the aid of Proposition 2, the proof of this result now follows as in \cite{mp}.

\end{document}